\documentclass[12pt]{amsart}
\usepackage{amsthm, amssymb, a4wide}
\usepackage[arrow,matrix]{xy}
\usepackage{pstricks,pst-node, pst-grad}
\usepackage{hyperref}
 \hypersetup{
    pdftitle=   {The number of orientable small covers over cubes},
    pdfauthor=  {Suyoung Choi},
    pdfsubject={Orientable small covers over cubes},
  pdfkeywords={Orientable small cover, acyclic digraph, real Bott manifold, Toric topology}
}

\theoremstyle{plain}
\newtheorem{theorem}{Theorem}[section]

\newtheorem{corollary}[theorem]{Corollary}
\theoremstyle{definition}
\newtheorem{remark}[theorem]{Remark}
\newcommand{\R}{\mathbb{R}}
\newcommand{\Z}{\mathbb{Z}}
\newcommand{\RP}{\mathbb{R}P}

\def\iput(#1)#2{\cnode[linewidth=0.7pt](#1){5pt}{#2}\rput(#1){{$#2$}}}

\title[Orientable small covers over cubes]{The number of orientable small covers over cubes}
\author[S.Choi]{Suyoung Choi}
\address{Department of Mathematics, Osaka City University, Sugimoto, Sumiyoshi-ku, Osaka 558-8585, Japan}
\email{choi@sci.osaka-cu.ac.jp}
\urladdr{http://math01.sci.osaka-cu.ac.jp/~choi}
\subjclass[2000]{Primary 37F20, 57S10; Secondary 57N99}
\keywords{Orientable small cover, acyclic digraph, real Bott manifold, Toric topology}

\begin{document}
\maketitle

\begin{abstract}
We count orientable small covers over cubes. We also get estimates for $O_n/R_n$, where $O_n$ is the number of orientable small covers and $R_n$ is the number of all small covers over an $n$-cube up to the Davis-Januszkiewicz equivalence.
\end{abstract}

\section{Introduction}
A \emph{small cover}, defined by
Davis and Januszkiewicz \cite{DJ}, is an $n$-dimensional
closed smooth manifold $M$ with an effective real
torus $(S^0)^n(=: T^n)$-action such that the action is
locally isomorphic to a standard $T^n$-action on $\R^n$
and the orbit space $M/T^n$ can be identified with a
simple combinatorial polytope. For instance, $\RP^n$
with a natural $T^n$-action is a small cover over an $n$-
simplex. In general, a real toric manifold, the set of
real points of a toric manifold, provides an example
of small covers. Hence, small covers can be seen as
a topological generalization of real toric manifolds in
algebraic geometry.

A small cover over a cube is known as a \emph{real Bott
manifold} which is obtained as iterated $\RP^1$ bundles
starting with a point, where each fibration is the projectivization of a Whitney sum of two real line bundles. These manifolds are well-studied in numerous
papers such as \cite{Ka-Ma-2009} and \cite{Masuda-2008}. The author also found a
strong relation between small covers and acyclic digraphs, and he calculated the number of them up to
several senses in \cite{Choi-2008}.

In the present paper, we restrict our attention
to the case of orientable small covers over a cube.
Thankfully, Nakayama and Nishimura \cite{Na-Ni-2005} found a
simple criterion for a small cover to be orientable.
Using this criterion, we establish the formula of the
number of orientable small covers over a cube and
show that the ratio $O_n / R_n$ is approximately $\frac{1.262}{2^n}$,
where $O_n$ is the number of orientable small covers
and $R_n$ is the number of small covers over an $n$-cube
up to the Davis-Januszkiewicz equivalence.

\section{Orientable small covers over cubes}
Let $P$ be an $n$-dimensional simple polytope with $m$ facets.
Two small covers $M_1$ and $M_2$ over $P$
are \emph{Davis-Januszkiewicz equivalent} (or simply, \emph{D-J
equivalent}) if there is a weak $T^n$-equivariant homeomorphism $f \colon M_1 \rightarrow M_2$ which makes the diagram
commute:
$$
\xymatrix{
M_1 \ar[rr]^{f} \ar[dr] & & \ar[ld] M_2 \\
 & P&
}.
$$

It is well-known by \cite{DJ}  that all small covers over
$P$ can be distinguished by the map $\lambda$ from the set of
facets of $P$ to $\Z_2^n=\{0,1\}^n$, called the \emph{characteristic function}, which satisfies the \emph{non-singularity condition}; $\{\lambda(F_{i_1}), \ldots,
 \lambda(F_{i_n})\}$ is a basis of $\Z_2^n$ whenever
the intersection $F_{i_1} \cap \cdots \cap F_{i_n}$ is non-empty, where $\{F_1, \ldots, F_m\}$ is the set of facets of $P$. Let $M_1,~M_2$ be two small covers over $P$ corresponding to characteristic functions $\lambda_1,~\lambda_2$, respectively. By \cite{DJ}, $M_1$ is D-J equivalent to $M_2$ if and only if there is an automorphism $\sigma \in \text{Aut}(\Z_2^n)$ such that $\lambda_1 = \sigma \circ \lambda_2$. Hence,
the D-J equivalence classes are independent of the
choice of basis for $\Z_2^n$. One may assign an $n \times m$ matrix $\Lambda$ to $\lambda$ by ordering the facets and choosing a basis for $\Z_2^n$ as the follow:
$$
    \Lambda = \left( \lambda(F_1) \cdots \lambda(F_m) \right).
$$
If we additionally assume that the first $n$ facets meet
at a vertex, by the non-singularity condition, we can
choose an appropriate basis of $\Z_2^n$ such that $\Lambda = (E_n | \Lambda_\ast)$, where $E_n$ is the identity matrix of size $n$
and $\Lambda_\ast$ is an $n \times (m-n)$ matrix. Hence, the D-J
equivalence classes of small covers over $P$ are classified by $\Lambda_\ast$.

Now, we consider the case where $P$ is an $n$-cube.
Note that $P$ has $2n$ facets. We order the facets of
$P$ satisfying $F_j \cap F_{n+j} = \emptyset$ for $1 \leq j \leq n$. Then the first $n$ facets meet at a vertex. Hence, for each $\lambda$, the corresponding matrix $\Lambda$ can be expressed as $\Lambda = (E_n | \Lambda_\ast)$, where $\Lambda_\ast$ is an $n\times n$ matrix. One can check that the non-singularity condition holds if and only if all of principal minors of $\Lambda_\ast$ are $1$. Therefore, there is a bijection between small covers over cubes up to the D-J equivalence and square $\Z_2$-matrices all of whose principal minors are $1$.

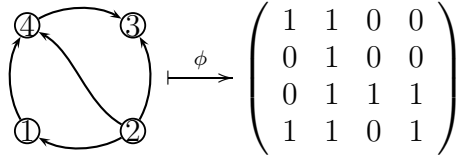
\begin{figure}
\psset{unit=10pt}
\[
\xymatrix{
\begin{pspicture}(0,0)(6,6)
\iput(1,1){1} \iput(5,1){2} \iput(5,5){3} \iput(1,5){4}
\nccurve[angleA=210,angleB=-30]{->}{2}{1}
\nccurve[angleA=60,angleB=300]{->}{2}{3}
\nccurve[angleA=150,angleB=330]{->}{2}{4}
\nccurve[angleA=120,angleB=240]{->}{1}{4}
\nccurve[angleA=30,angleB=150]{->}{4}{3}
\end{pspicture}
\ar@{|->}[r]^-{\phi} &
{\left(
  \begin{array}{cccc}
    1 & 1 & 0 & 0 \\
    0 & 1 & 0 & 0 \\
    0 & 1 & 1 & 1 \\
    1 & 1 & 0 & 1 \\
  \end{array}
\right)}
}
\]
  \caption{A bijection $\phi$}\label{fig:bijection}
\end{figure}

Let $M(n)$ be the set of square $\Z_2$-matrices of
size $n$ all of whose principal minors are $1$ and let $\mathcal{G}_n$
be the set of acyclic digraphs with labelled $n$ vertices.
By \cite{Choi-2008}, we have a bijection $\phi : \mathcal{G}_n
\rightarrow M(n)$ by
$$
    \phi : G \mapsto A(G)^t+E_n,
$$ where $A(G)^t$ is the transpose matrix of the vertex adjacency matrix of $G$ (see Figure~\ref{fig:bijection}).

\begin{remark} In the classical theory of real
Bott manifolds, the representative matrix of real
Bott manifold is the transpose matrix of its characteristic function matrix $\Lambda_\ast$. This is a reason why
we use $A(G)^t$ instead of $A(G)$ in the definition of $\phi$.
\end{remark}

On the other hand, we have a nice orientability condition for small covers due to Nakayama and Nishimura in \cite{Na-Ni-2005}.

\begin{theorem}[Nakayama and Nishimura \cite{Na-Ni-2005}] \label{theorem:orientablility condition}
Let $P$ be an $n$-dimensional simple polytope with $m$ facets and let $M$ be a small cover over $P$ with $\Lambda$. Then $M$ is
orientable if and only if the sum of entries of the $i$-th column of $\Lambda$ is odd for all $i = 1, \ldots , m$.
\end{theorem}

\begin{corollary} \label{lemma:orientable small cover}
The number of orientable
small covers over an $n$-cube up to D-J equivalence is
equal to the number of acyclic digraphs with labelled
$n$ vertices all of whose vertices have even out-degrees.
\end{corollary}
\begin{proof}

    Let $G$ be a digraph and $A(G)$ its vertex adjacency matrix. Then the sum of entries of the $i$-th row of $A(G)$ means the out-degree of the $i$-th vertex of $G$ (see Appendix). Let $M$ be a orientable small cover over an $n$-cube corresponding to $\Lambda_\ast$. Since $\Lambda_\ast \in M(n)$, the transpose $\Lambda_\ast^t$ of $\Lambda_\ast$ is also in $M(n)$. Note that the sum of entries of each row of $\Lambda_\ast^t - E_n$ is even by Theorem \ref{theorem:orientablility condition}, and hence, every vertex of $\phi^{-1}(\Lambda_\ast)$ has an even out-degree. Since the D-J equivalence classes are classified by $\Lambda_\ast$ and $\phi$ is a bijection, we prove the corollary.
\end{proof}

\section{The number of orientable small covers}

Let $R_n$ be the number of acyclic digraphs with labelled $n$
vertices. The following is the recursive formula for $R_n$ due to
R.~W.~Robinson in \cite{MR0276143}.
\begin{equation*}\label{equation:R_n}
    R_n = \sum_{k=1}^n (-1)^{k+1} {n \choose k} 2^{k(n-k)} R_{n-k}.
\end{equation*}

Let $\mathcal{O}_n \subset \mathcal{G}_n$ be the set of acyclic digraphs all of whose vertices have even out-degrees and let $O_n$ be the
cardinality of $\mathcal{O}_n$ (we use the alphabet `O' instead of `E' although they have only `even' out-degree vertices, because the `O' is the abbreviation of the word `\emph{Orientable}').

\begin{theorem}\label{theorem:O_n}
Let $R_k$ be the number of acyclic digraphs with labelled $k$
vertices. Then,
$$
        O_n = \sum_{k=1}^n (-1)^{k+1} {n \choose k} 2^{(k-1)(n-k)} R_{n-k}.
$$
\end{theorem}
\begin{proof}
We count matrices in $M(n)$ all of whose the sum of entries of each
column are odd. Let us denote the sum of entries of the $i$-th
column of an $n\times n$ matrix $A$ by $c_i(A)$. Since an acyclic
digraph always has a vertex of out-degree $0$, there is at least one $i$ such that
$c_i(A)=1$ for each $A \in M(n)$. Assume $c_{i_1}(A)=\cdots =c_{i_k}(A)=1$,
where $k \geq 1$. Since all principal minors of $A$ are $1$, the
diagonal entries of $A$ are all $1$. Thus, by a replacement of labels, we may
assume that $A$ is of the following form:
\begin{equation}\label{equation:I-E principle}
    \left(
         \begin{array}{cc}
            E_k & S \\
            0 & T \\
         \end{array}
    \right),
\end{equation} where $E_k$ is the identity matrix of size $k$, $T$ is an $(n-k)\times (n-k)$-matrix and $S$ is a $k \times (n-k)$-matrix. Note that $A \in M(n)$ if and only if $T \in M(n-k)$. Thus we may control only one row of $S$ for making all $c_i(A)$'s are odd. This implies the number of $A$'s of the form (\ref{equation:I-E principle}) whose $c_i(A)$'s are odd for all $i$ is $2^{(k-1)(n-k)} R_{n-k}$. To avoid counting repeatedly, we apply the Principle of Inclusion-Exclusion and we get the formula for $O_n$.
\end{proof}

Here are a few values of $R_n$ and $O_n$.

\begin{center}
\begin{tabular}{c|c c c c c c c }
  \hline
  $n$ & 1 & 2 & 3 & 4 & 5 & 6 & 7 \\ \hline
  $R_n$ &1 &3& 25 &543 &29,281& 3,781,503 &1,138,779,265\\
  $O_n$ & 1 & 1 & 4 & 43 & 1,156 & 74,581 & 11,226,874 \\
  \hline
\end{tabular}
\end{center}

Let us consider the \emph{chromatic generating functions} of $R_n$ and $O_n$, namely, we set
$$
    R(x) = \sum_{n=0}^\infty R_n \frac{x^n}{n! 2^{n \choose 2}} \hbox{ and } O(x) = \sum_{n=0}^\infty O_n \frac{x^n}{n! 2^{n \choose 2}}.
$$
\begin{corollary}
    Let $F(x) = \sum_{n=0}^\infty \frac{x^n}{n! 2^{n \choose 2}}$. Then
$$
    O(x)= \frac{1 - F(-x)}{F(-\frac{x}{2})}.
$$
\end{corollary}
\begin{proof}
Let us consider chromatic generating functions $A(x),~B(x)$ and $C(x)$ with respect to the sequences $A_n,~B_n$ and $C_n$, respectively. Note that if $C(x)=A(x)B(x)$, then $C_n = \sum_{k=0}^n A_k B_{n-k} {n \choose k} 2^{k(n-k)}$. Thus, we have $F(-x)R(x)=1$ (see \cite{MR0317988}) and
$$
    R\left( \frac{x}{2} \right)F(-x) + O(x) = \sum_{n=0}^\infty \frac{R_n}{2^n} \frac{x^n}{n! 2^{n \choose 2}} = R\left( \frac{x}{2} \right).
$$ Hence we have $O(x) = F\left(- \frac{x}{2} \right)^{-1}(1-F(-x))$.
\end{proof}

Let $G(x) = \frac{F(\frac{x}{2})}{1-F(x)}$. We obtain estimates
for $O_n$ by analyzing the behavior of the function $G(x)$. Since
$F(x)$ has an isolated zero $\alpha \approx -1.488$ (see \cite[Section 2]{Stanley-1997}), $G(x)$ has an
isolated zero $2\alpha$. Hence, standard techniques provide the asymptotic formula
$$ G(x) \sim G'(2\alpha)(x-2\alpha).$$
Hence we have
$$
    O(x) = \frac{1}{G(-x)} \sim \frac{1}{G'(2\alpha)(-x-2\alpha)}.
$$
Note that $F'(x) = F(\frac{x}{2})$. Therefore, the following asymptotic formula
$$
   O(x) \sim - \frac{1-F(2\alpha)}{\alpha F(\frac{\alpha}{2})} \sum_{n=0}^\infty
    \left(-\frac{x}{2\alpha}\right)^n
$$
immediately follows two facts $\frac{1}{G'(2\alpha)} = \frac{2(1-F(2\alpha))}{F'(\alpha)}$ and $\frac{1}{-x-2\alpha} = \frac{1}{-2\alpha} \sum_{n=0}^\infty (- \frac{x}{2\alpha})^n$.

Therefore $O_n \sim K 2^{\binom{n}{2}}n! \left(-
\frac{1}{2\alpha}\right)^n$, where $K = -
\frac{1-F(2\alpha)}{\alpha F\left(\frac{\alpha}{2}\right)}
\approx 2.197$.

\begin{corollary} \label{cor:ratio}
We have estimates for the orientable small covers ratio as
$$
    \frac{O_n}{R_n} \sim \frac{K}{C 2^n},
$$ where $\frac{K}{C} \approx 1.262$.
\end{corollary}%
\begin{proof}
Since $R(x)F(-x)=1$ and $F(x)$ has an isolated zero $\alpha$,
we have $R(x) =\frac{1}{F(-x)} \sim \frac{1}{F'(\alpha)(-x-\alpha)} = \frac{1}{-\alpha F(\frac{\alpha}{2})} \sum_{n=0}^\infty (- \frac{x}{\alpha})^n$.
Hence, we
have $R_n \sim C 2^{\binom{n}{2}}n!\left(-
\frac{1}{\alpha}\right)^n$, where $C =
-\frac{1}{\alpha F\left(\frac{\alpha}{2}\right)}\approx 1.739$. Therefore
$\frac{O_n}{R_n} \sim \frac{K}{C 2^n} \approx \frac{1.262}{2^n}$.
\end{proof}

\section*{Appendix. Graph theory terminology}
We review the terminology in graph theory, following
\cite{Stanley-1997}. A \emph{directed graph} or \emph{digraph} $G$
is a triple $(V,E,\varphi)$, where $V=\{v_1, \ldots, v_n\}$ is a set
of vertices, $E$ is a set of directed edges, and
$\varphi$ is a map from $E$ to $V \times V$.
If $\varphi(e)=(u,v)$, then $e$ is
called an \emph{edge from $u$ to $v$} with the initial vertex $u$
and the final vertex $v$. If $u=v$ then $e$ is called a
\emph{loop}. If $\varphi$ is injective and has no loops, then $G$ is said to be \emph{simple}. In this case, we denote $e$ by $(u,v)$ for
simplicity and represent $G$ by $(V,E)$. \emph{Throughout
this paper, every graph is simple.} A \emph{walk} of length $k$
from vertex $u$ to $v$ is a sequence $v_0, v_1, \ldots, v_k$ such
that $v_0 = u$ and $v_k = v$, where $ (v_i, v_{i+1}) \in E$ for all
$i=0, \ldots, k-1$. If all the $v_i$'s are distinct except for $v_o
= v_k$, then the walk is called a \emph{cycle}. $G$ is
\emph{acyclic} if there is no cycle of any length in $G$. The \emph{out-degree}
of a vertex $v$ is the number of edges of $G$ with the
initial vertex $v$. Similarly the in-degree of $v$ is the
number of edges of $G$ with the final vertex $v$.

All digraphs can be represented by matrices.
Define an $n \times n$ matrix $A(G) = (A_{ij})$ by
 $$
    A_{ij} =\left\{
              \begin{array}{ll}
                1, & \hbox{if $(v_i, v_j) \in E$;} \\
                0, & \hbox{otherwise.}
              \end{array}
            \right.
$$ The matrix $A(G)$ is called the \emph{vertex adjacency matrix} of $G$. We remark that the sum of entries of the $i$-th column of $A$ is equal to the in-degree of $v_i$ and the sum of entries of the $j$-th row of $A$ is equal to the out-degree of $v_j$.

\section*{Acknowledgements}
The research of the author was carried out with the support of the Brain Korea 21 project, KAIST, and the Japanese Society for the Promotion of Sciences (JSPS grant no. P09023).
The author is deeply grateful to Professor Mikiya
Masuda for bringing a ratio question (Corollary~\ref{cor:ratio}) to my attention when they visited at Fudan University in
January 2008.

\end{document}